\documentclass[final,leqno]{siamltex}

\usepackage{amssymb,amsfonts,amsmath}                        % math symbols
\usepackage{bm}
\usepackage{graphicx}
\usepackage{caption}
\usepackage{subcaption}
\usepackage[section]{placeins}
\usepackage{booktabs}
\title{Numerics for Hyperbolic Conservation Laws with Help from the physical Entropy}

\author{Carl Philipp Zinner\thanks{ETH Z\"urich, Department of Materials, Polymer Physics, HCI H 538, Vladimir-Prelog-Weg 1-5/10, CH-8093 Z\"urich, Switzerland ({\tt carl.zinner@mat.ethz.ch; http://www.polyphys.mat.ethz.ch/}).}
\and Hans Christian \"Ottinger\thanks{ETH Z\"urich, Department of Materials, Polymer Physics, HCI H 543, Vladimir-Prelog-Weg 1-5/10, CH-8093 Z\"urich, Switzerland  ({\tt hco@mat.ethz.ch; http://www.polyphys.mat.ethz.ch/}).}}

\begin{document}

\maketitle

\begin{abstract}
Stable numerical simulations for a hyperbolic system of conservation laws of relaxation type but not in divergence form are obtained by incorporating the physical entropy into the simulations. The entropy balance is utilized as an additional equation to eliminate the numerically critical terms with simple substitutions. The method has potential for a wider applicability than the particular example presented here.
\end{abstract}



\pagestyle{myheadings}
\thispagestyle{plain}
\markboth{C.P. ZINNER and H.C.  \"OTTINGER}{NUMERICS FOR CONSERVATION LAWS WITH ENTROPY}

\section{Introduction}
\label{sec:intro}
The advantages of a thermodynamically admissible system with a positive entropy production are obvious from a physical point of view. For simulations, the entropy balance is often regarded as redundant. Motivated by the success of the Entropic Lattice Boltzmann Method \cite{karlin2006elements},\cite{karlin2007elements}, where the entropy acts as a stabilizing Lyapunov function, we seek to incorporate the entropy balance in a shock tube study of moment equations for rarefied gas dynamics. The modified system remains well posed by removing terms from the original system. This elimination enables us to rewrite the system of equations in divergence form and use a tailored numerical scheme for hyperbolic conservation laws. The stability of the resulting simulations encourages us to use this straightforward recipe for other entropic models with numerically challenging properties.

The paper is structured as follows. We start by introducing the nomenclature and mathematical generalization of the system that we intend to study and give a brief overview of the available numerical methods in section \ref{sec:hyp}. Section \ref{sec:hco} outlines the physical model of interest. In section \ref{sec:tube} we discuss the intended simulation setup for a shock tube simulation. Section \ref{sec:num} describes the family of numerical schemes that we are building on and a naive extension for the numerically problematic terms in our physical model. A more sophisticated method using the entropy to overcome the gap in the numerical literature is derived in section \ref{sec:sol}. Finally, section \ref{sec:results} presents the simulation results and gives a comparison between the naive and entropic extension.

\section{Hyperbolic conservation laws}
\label{sec:hyp}
Consider a system of first order partial differential equations that constitute conservation laws of production type in one spatial dimension,
\begin{equation}\label{eq:div}
\frac{\partial u}{\partial t} + \frac{\partial f(u)}{\partial x} = \frac{1}{\tau}g(u)
\end{equation}
where $u \in \mathbb{R}^m$ is the vector of \textit{conserved} variables, $f(u) : \mathbb{R}^m \rightarrow \mathbb{R}^m$ is the flux vector, $g(u): \mathbb{R}^m \rightarrow \mathbb{R}^m$ is a production term and $\tau$ a relaxation constant. Such a system is called \textit{hyperbolic} in $t$ at $u$ if the eigenvalues of the quasilinear Jacobian $\frac{\partial f(u)}{\partial u}$ are real. In physical terms, the hyperbolicity condition means that the characteristic speeds of the quasilinear equations are finite. This is in contrast to parabolic partial differential equations with imaginary eigenvalues, such as the Navier-Stokes equations, which can be seen as a regularization of the hyperbolic Euler equations. The finite speeds pose a numerical challenge for simulations, since a naive approximation of the derivatives does not necessarily reproduce the information propagation accurately. Problems arise in particular for non-smooth solutions such as shock structures in gas dynamics. These are discontinuities that travel at supersonic speeds, such as a pressure shock at the wing edges of a supersonic jet.

Many numerical schemes have been developed for the homogenous case $g(u)=0$ in the context of the Euler equations of gas dynamics. A comprehensive introduction to the numerical methods available is given in \cite{leveque1992numerical}.  The majority of these schemes make use of the characteristic structure of the equations by solving Riemann problems on the discrete level, since the eigenvalue structure of the Euler equations is well understood. For some systems the eigenvalues cannot be obtained analytically. One method able of handling this situation is the Nessyahu-Tadmor \cite{nessyahu1990non} central scheme. This method has been extendend for the non-homogeneous case of production terms \cite{liotta2000central}, where a family of explicit and implicit schemes for stiff relaxation are proposed. This scheme has been successfully implemented for various moment equations in gas dynamics and is the starting point of our numerical approach.

The mathematics of hyperbolic conservation laws are no longer well understood, as soon as the equations cannot be written in divergence form as above in (\ref{eq:div}). Naturally this carries over to the numerical methods and there is a lack of schemes for this type of hyperbolic systems. We present such a system in the following section.

\section{Thermodynamically admissible 13 moment equations from the Boltzmann equation}
\label{sec:hco}
The physical model we are studying is a set of moment equations in three space dimensions derived from the Boltzmann equation. Its purpose is to describe rarified monatomic gas dynamics and gas flows far from equilibrium where the constitutive equations of Navier-Stokes-Fourier no longer hold (\cite{paperhco}, \cite{hcocomment} and \cite{ottinger2010ottinger}). These rarefied flow domains are best characterised by a dimensionless Knudsen number ${\rm{Kn}}=\frac{l_{\rm{mfp}}}{L}$ of order unity, which refers to the ratio of the particle mean free path $l_\text{mfp}$ over the characteristic length of the problem $L$. A classical example of such a flow situation is the initial re-entry phase of a space vehicle.

Moment equations reach a more appropriate level of description by extending the five fields of hydrodynamics by a second moment tensor $\boldsymbol{\pi}$ and third moment vector $\boldsymbol{q}$, which are related to the pressure tensor and heat flux, to a total of thirteen fields.
In this section we present and discuss the fully three dimensional system of equations before reducing the system to the one dimensional simulation setup. For the density $\rho$ and momentum $\boldsymbol{M}$ we have the common mass and momentum balance
\begin{equation}\label{eq:density}
\frac{\partial \rho}{\partial t} + \frac{\partial M_k}{\partial x_k} = 0,
\end{equation}
and
\begin{equation}
 \frac{\partial M_i}{\partial t} + \frac{\partial }{\partial x_k} \left(\frac{ M_i M_k}{\rho} + \frac{\rho}{m^2} \pi_{ik}\right) = 0,
\end{equation}
where we have used Einstein's summation convention.
The additional second moment tensor $\boldsymbol{\pi}$ and skewed third moment vector $\boldsymbol{q}$ are governed (\cite{paperhco} and \cite{ottinger2010ottinger}) by the evolution equations
\begin{multline}\label{eq:pi}
 \frac{\partial \pi_{ij}}{\partial t} = -v_k \frac{\partial \pi_{ij} }{\partial x_k} - \frac{\partial v_i}{\partial  x_k} \pi_{kj}  - \pi_{ik}\frac{\partial v_j}{\partial x_k}\\
 - \frac{F}{m} \frac{2 b}{1+ 2 \varphi b} \left(q_i\frac{\partial \pi_{kk}}{\partial x_j} + \frac{\partial \pi_{kk}}{\partial x_i} q_j\right)\\
 - 2\frac{\pi_{ij}}{\rho}\frac{\partial }{\partial x_k} \left(\frac{F}{m} \frac{2 b}{1+ 2 \varphi b} \rho q_k\right)- \frac{1}{\tau} \left(\pi_{ij} - \frac{\pi_{kk}}{3} \delta_{ij}\right),
\end{multline}
and
\begin{equation}\label{eq:q}
 \frac{\partial q_i}{\partial t} = -v_k \frac{\partial q_i}{\partial x_k}  - \frac{\partial v_i}{\partial x_k} q_k -\frac{F}{m} \frac{\partial \pi_{kk}}{\partial x_i}\\
 -\frac{q_i}{\rho}\frac{\partial }{\partial x_k} \left(\frac{F}{m} \frac{2 b}{1+ 2 \varphi b} \rho q_k\right) - \frac{1}{2 \tau} D_{ik} q_k,
\end{equation}
where $\boldsymbol{v} = \boldsymbol{M}/\rho$ is the velocity, $m$ is the particle mass and $\varphi = \boldsymbol{q} \cdot \boldsymbol{\pi}^{-1} \cdot \boldsymbol{q}$ is an often occurring scalar indicating the skewness away from equilibrium. The tensor $\boldsymbol{D}= \boldsymbol{\bar{D}} + 1 -(\rm{tr}\boldsymbol{\pi}/3)\boldsymbol{\pi}^{-1}$ contains the isotropic dimensionless tensor $\boldsymbol{\bar{D}}$ which is determined by the two-particle interaction potential. The temperature (and pressure) is contained in the trace of the second moment tensor. The dimensionless constants $F$ and $b$ are related to the closure, and $\tau$ is a relaxation time.
In comparison to the well known 13 moment equations of Grad \cite{grad}, this system has the advantage of being thermodynamically admissible and equipped with an entropy
\begin{equation}\label{eq:entropy}
 S = \frac{k_\text{B}}{m} \int_V \left\{\frac{1}{2} \ln \left[ \left(\frac{2 \pi}{h^2}\right)^3 \frac{m^2}{\rho^2} \det \boldsymbol{\pi} \right] + \frac{5}{2} + \bar{S}(\varphi)\right\} \rho d^3 x.
\end{equation}
Here we have introduced the skewness contribution $\bar{S}(\varphi)=-b\varphi$ and the Boltzmann and Planck constants $k_\text{B}$ and $h$. The evolution of the entropy density $s$ is governed by
\begin{equation}
 \frac{\partial s}{\partial t} + \frac{\partial }{\partial x_k} \left(v_k s + \frac{3 k_\text{B}}{m^2} \frac{F}{m} \frac{2 b}{1+ 2 \varphi b} \rho q_k\right) =
 \frac{n k_\text{B}}{6 \tau}(\pi_{kk}\pi_{kk}^{-1} -9) + \frac{n k_\text{B}b}{\tau} q_i \pi^{-1}_{ij}\bar{D}_{jk}q_k,
\end{equation}
where $n$ is the particle number density $n=\rho/m$.

A detailed study of the physical properties of this system, including linear and nonlinear stability analysis, hyperbolicity, sound waves, shock structures, boundary conditions and channel flows is subject to forthcoming papers. In the following we discuss the numerical properties of the system and only give a brief description of one of the applications: the shock tube.

\section{Shock tube}
\label{sec:tube}
The most prominent test case for the physical accuracy of gas dynamics models and hyperbolic numerical schemes is the one-dimensional shock tube problem. In this application we consider a long tube separated into two sections by a membrane. The gases in the two chambers are in equilibrium, but differ in pressure and density, where the pressure difference can be up to two
orders of magnitude. After the membrane is burst, one can observe a shock wave and a contact discontinuity traveling at supersonic speeds into the low pressure section of the tube, and a rarefaction wave traveling into the high pressure section.
Apart from its interesting physical properties, this application is the standard benchmark for the robustness of numerical methods due to the inherent non-smoothness \cite{leveque1992numerical}. Mathematically, the discontinuous initial conditions constitute a Riemann problem, of which the analysis is a vital element of the theory of hyperbolic equations and their numerical methods.
Additionally, the geometry in Cartesian coordinates is interpreted as essentially one-dimensional and therefore no boundary conditions are needed.

The relevant spatial dimension $x_1$ describes the position along the tube and shows the only dynamical variation. The tube's width and height (or circular cross section) and the size of the membrane are infinite and therefore the problem is independent of $x_2$ and $x_3$. This leads to the state variables density $\rho$, momentum $M_1$ in direction $x_1$ and the heat flux related $q_1$. The radial symmetry now imposes a second moment tensor of
\begin{equation}
\boldsymbol{\pi} = \left(
\begin{array}{ccc}
\pi_{11} & 0 & 0 \\
0 & \pi_{22} & 0 \\
0 & 0 & \pi_{22}
\end{array} \right).
\end{equation}
The next step is to find a diffeomorphic transformation of variables in order to rewrite the field equations (\ref{eq:density} - \ref{eq:q}) in the desired balance form (\ref{eq:div}). Naturally, the guiding principle is conservation, and density and momentum are already conserved. The second moment tensor (\ref{eq:pi}) itself is not conserved but its trace can be incorporated in a conserved equation for the mechanical energy
\begin{equation}
 E= \int_V \left(\frac{\boldsymbol{M}^2}{2\rho} + \frac{\rho}{2m^2}\text{tr}\boldsymbol{\pi}\right)d^3x.
\end{equation}
Additionally we take the second diagonal term of the non-peculiar second moment tensor (including a momentum contribution) which includes the deviatoric stress. Finally, for the remaining variable we take $q_1$ as it is. This transformation yields the following system of equations: in the notation of (\ref{eq:div}) the conversion between the simulation variables $\boldsymbol{u}$ and the physical variables $\boldsymbol{w}$ are given in Table \ref{table:vars}.
\begin{table*}[!hp]\centering
\begin{tabular}{@{}lrclr@{}}\toprule
\multicolumn{2}{c}{simulation} & \phantom{abcdefgh}& \multicolumn{2}{c}{$\textrm{physical}$}\\
\cmidrule{1-2} \cmidrule{4-5}
$\boldsymbol{u}$ & $\boldsymbol{w}$ && $\boldsymbol{w}$ & $\boldsymbol{u}$\\
\midrule
$u_1$ & $\rho$ && $\rho$ & $u_1$\\
$u_2$ & $M_1$ && $M_1$ & $u_2$ \\
 \phantom{ab} & \phantom{ab} && $v_1$ & $\frac{u_2}{u_1}$\\
$u_3$ & $ \frac{1}{2} \rho v_1^2 + \frac{1}{2m^2} \rho (\pi_{11}+2\pi_{22})$ && $\pi_{11}$ & $\frac{m^2(-u_2^2+ 2 u_1 u_3-4u_1 u_4)}{u_1^2}$\\
$u_4$ & $\frac{1}{2m^2} \rho \pi_{22}$ && $\pi_{22}$ & $\frac{2 m^2 u_4}{u_1}$\\
$u_5$ & $q_1$ && $q_1$ & $u_5$\\
\bottomrule
\end{tabular} \caption{Variable conversion}\label{table:vars} \end{table*}

The flux vector $f(u)$ is given by the components
\begin{equation}
\begin{array}{l}
f_1=M_1 \\
[1.5pt]\\
f_2=v_1 M_1 + \dfrac{\rho}{m^2} \pi_{11} \\
[1.5pt]\\
f_3=\dfrac{1}{2} \rho v_1^3 +  \dfrac{1}{2} \dfrac{1}{ m^2} \rho v_1  (3 \pi_{11} +2 \pi_{22})+ \dfrac{1}{m^2} \dfrac{F}{m} \dfrac{2 b}{1+ 2 \varphi b} \rho \pi_{11}q_1 + 2 \dfrac{1}{m^2} \dfrac{F}{m} \dfrac{2 b}{1+ 2 \varphi b} \rho \pi_{22}q_1 \\
[1.5pt]\\
f_4=\dfrac{1}{2 m^2} \rho \pi_{22} v_1  \\
[1.5pt]\\
f_5=v_1 q_1 + \dfrac{F}{m}(\pi_{11}+2\pi_{22})
\end{array}
\end{equation}
and the production vector is given by
\begin{equation}
\begin{array}{l}
g_1 = g_2 = g_3 =0 \\
[1.5pt]\\
g_4 = -\dfrac{1}{2 m^2} \rho \dfrac{1}{3}\dfrac{1}{\tau}(\pi_{22}-\pi_{11})  \\
[1.5pt]\\
g_5 = -\dfrac{1}{2 \tau}D_{11}q_1
\end{array}
\end{equation}
where we have included the relaxation parameter $\tau$.

The variable transformation has brought us as close as possible to a clean divergence form, but leaves remaining terms on the right hand side that are neither of the production type, nor the flux-divergence type:
\begin{equation}\label{eq:remainder}
\begin{array}{l}
h_1 = h_2= h_3 =0\\
[1.5pt]\\
h_4 =- \dfrac{1}{m^2} \pi_{22} \dfrac{\partial }{\partial x_1} \left( \dfrac{F}{m} \dfrac{2 b}{1+ 2 \varphi b} \rho q_1 \right) \\
[1.5pt]\\
h_5 = -\dfrac{q_1}{\rho} \dfrac{\partial }{\partial x_1} \left( \dfrac{F}{m} \dfrac{2 b}{1+ 2 \varphi b} \rho q_1 \right).
\end{array}
\end{equation}
The derivatives in (\ref{eq:remainder}) pose a numerical challenge and we propose a simple method to overcome them.

\section{Numerical method}
\label{sec:num}
Before we find a solution for the additional terms of our system, we will address the other numerical properties in order to select the appropriate scheme. These are:
\begin{enumerate}
\item Hyperbolicity.
\item An eigenstructure that is only known from numerical evaluation.
\item A relaxation constant $\tau$ that can vary over a wide range of regimes of rarefaction and hence considerable stiffness.
\end{enumerate}
The scheme addressing the first two points is the central Nessyahu-Tadmore scheme \cite{nessyahu1990non}. It is a predictor-corrector method of second degree accuracy in time and space and works with and without a detailed knowledge of the eigensystem. The extension for stiff production terms is a family of implicit schemes by \cite{liotta2000central}.

Let us define the \textit{minmod} function
\begin{equation}
\text{mm}(x,y) = \left\{
\begin{array}{rl}
\operatorname{sgn} (x) \min (|x|, |y|) & \text{if } \operatorname{sgn} (x) = \operatorname{sgn} (y),\\
0 & \text{otherwise}
\end{array} \right.
\end{equation}
which is required for the discrete derivatives $\partial u_j/\partial x \approx u_j'/\Delta x$ at spatial grid nodes $j$ and cell size $\Delta x$:
\begin{align}
 f_j'  & = \text{mm}(f_{j+1} - f_j -\frac{1}{2} D_{j+\frac{1}{2}}f, f_{j} - f_{j-1} +\frac{1}{2} D_{j-\frac{1}{2}}f),\\
D_{j+\frac{1}{2}}f  &= \text{mm}(f_{j+2} - 2 f_{j+1} +f_{j}, f_{j+1} - 2 f_{j} + f_{j-1}).
\end{align}
These limiters are proposed by \cite{liotta2000central} as giving the best results.
The full time step integrates the cell $[x_j,x_{j+1}] \times [t^n,t^{n+1}]$. The integration consists of two predictor steps at time $\frac{1}{3} \Delta t$ and $\frac{1}{2} \Delta t$ for $u^{n+1/3}_j$ and $u^{n+1/2}_j$ at the cell edges, where $n$ denotes the time step index. The predictors are implicit in time, with the unknown variables appearing in the production $g$, and can be written in vector notation as
\begin{align}\label{eq:pred}
 u^{n+1/3}_j & = u_j^n + \frac{\Delta t}{3} \left( g(u^{n+1/3}_j)- \frac{f'_j}{\Delta x}\right),\\
 u^{n+1/2}_j & = u_j^n + \frac{\Delta t}{2} \left( g(u^{n+1/2}_j)- \frac{f'_j}{\Delta x}\right).
\end{align}
The corrector
\begin{align}
 u^{n+1}_{j+1/2} & = \frac{1}{2} (u_j^n + u_j^{n+1}) + \frac{1}{8} (u_j' - u_{j+1}')\\
 & -\lambda \left( f(u^{n+1/2}_{j+1}) - f(u^{n+1/2}_{j})\right)\\
  & + \Delta t \left( \frac{3}{8} g(u^{n+1/3}_j) + \frac{3}{8} g(u^{n+1/3}_{j+1}) + \frac{1}{4} g(u^{n+1}_{j+1/2}) \right)
\end{align}
evaluates the mid-cell values $x_{j+1/2}$ on a spatially staggered grid with $\lambda = \Delta t / \Delta x$, and is also time implicit with the unknows again appearing in the production term. After a full time step, the original grid is reconstructed by average interpolation. A continuous staggering results in oscillations depending on the size of $\Delta x$.

As presented, we have the ideal scheme for a system in divergence form. The remaining issue $h(u)$ could be dealt with in a naive manner: a second order central discretisation of the derivative appearing in (\ref{eq:remainder}) and subsequent treatment as an additional production term. The drawbacks apart from the obvious hyperbolicity issues are the loss of the implicitness and no reliable knowledge of the method's region of stability. In addition a few shock-structure related issues arise, which will be highlighted in the results section (\ref{sec:results}).

\section{Improved numerical implementation with the entropy balance}
\label{sec:sol}
Since our physical model has an entropy, we would like to take advantage of that. The terms appearing in $h(u)$, which are re-occurring, contain the expression
\begin{equation}
J_s = \frac{\partial }{\partial x_1} \left( \frac{F}{m} \frac{2 b}{1+ 2 \varphi b} \rho q_1 \right)
\end{equation}
which is the non-convective entropy flux $j_s$ without the pre-factor $\frac{3 k_\text{B}}{m^2}$. This allows us to integrate the entropy balance into the simulations whilst keeping a closed system of equations. We proceed in the following way: replace $J_s$, making $h(u)=0$ and rewriting $g(u,J_s)$ as
\begin{equation}\label{eq:g}
\begin{array}{c}
g_1 = g_2 = g_3 = 0 \\
[1.5pt]\\
g_4 = -\dfrac{1}{2 m^2} \rho \dfrac{1}{3}\dfrac{1}{\tau}(\pi_{22}-\pi_{11}) - \dfrac{1}{m^2} \pi_{22} J_s \\
[1.5pt]\\
g_5 =-\dfrac{1}{2 \tau}D_{11}q_1 -\dfrac{q_1}{\rho} J_s.
\end{array}
\end{equation}
We have declared $J_s$ as an additional unknown to the system and therefore can add the entropy evolution to the system, since it is no longer redundant and closes the extended system. Some care has to be taken at this step, since although we have the unknowns $u_s=(u_1, u_2, u_3, u_4, u_5, J_s)$ now, the system is not simply the same as (\ref{eq:div}) with system size $m+1$. The newly formed variable $J_s$ only appears in the production terms (\ref{eq:g}) and (\ref{eq:js}) and we therefore do not have a time evolution equation for it, but rather an algebraic relation. The coupling to the other variables is constructed from the entropy balance by expressing the entropy density $s$ in terms of the other five state variables with (\ref{eq:entropy}). It takes the form
\begin{equation}\label{eq:js}
 \frac{\partial s(u)}{\partial t} +\frac{\partial f_s(u)}{\partial x_1} = g_s(u) -\frac{3 k_\text{B}}{m^2} J_s,
\end{equation}
with a convective entropy flux
\begin{equation}
f_s(u) = v_1 s(u)
\end{equation}
and entropy production
\begin{equation}
g_s(u) = \frac{nk_\text{B}}{6 \tau} \left[\left(\pi_{11}+2\pi_{22}\right)\left(\frac{1}{\pi_{11}} +\frac{2}{\pi_{22}}\right) -9\right] + \frac{n k_\text{B}b}{\tau} q_1 \frac{1}{\pi_{11}}\bar{D}_{11}q_1.
\end{equation}
In its continuous form, the substitution amounts to expressing the non-convective entropy flux in a more complicated manner with the entropy balance. The fundamental information of the system has not changed. On a discretized level, this is a completely different story. What we have gained, is that the system is in divergence form and we can rely on the numerical recipe. The entropy balance turns into a discrete formula, and shows us how the problematic terms on the right hand side have to be treated in order to be consistent with the stable scheme. Although we lack the mathematical foundation to handle these terms in a hyperbolic system, we have made use of physical principles, a non-negative H-theorem, to enforce the stability of the simulations.
To illustrate the principle we write the modified first predictor step (\ref{eq:pred}).
\begin{align}
\frac{3 \left( u^{n+1/3}_j - u^{n}_j \right)}{\Delta t} + \frac{f'_j}{\Delta x} & = g(u^{n+1/3}_j,{J_s}^{n+1/3}_j),\\
\frac{ 3 \left( s(u^{n+1/3}_j) -s(u^{n}_j)\right)}{\Delta t} +\frac{{f'_s}_j}{\Delta x} & = g_s( u^{n+1/3}_j) -\frac{3 k_\text{B}}{m^2} {J_s}^{n+1/3}_j
\end{align}
This is still a system of six equations and six unknowns, namely $u^{n+1/3}_j$ and ${J_s}^{n+1/3}_j$, which can be solved at each time step. The new variable $J_s$ is still implicitly defined. It occurs in the production terms of the state variables $u$, explicitly on the left hand side of the additional equation and implicitly on the right hand side by way of the entropy.

The solutions $u^{n+1/3}_s=(u^{n+1/3}_1, u^{n+1/3}_2, u^{n+1/3}_3, u^{n+1/3}_4, u^{n+1/3}_5, J^{n+1/3}_s)$ (and $n+1/2$ as well as $n+1$ respectively) of this discrete system can be found by a root-finding algorithm for each predictor and corrector step. This can be achieved in two approaches. The slightly simpler approach is to solve all six unknowns by the root finder. The resulting program is easier to implement and modify, but the solver and starting values need some attention. Alternatively one can reformulate the predictor and corrector to have explicit formulas for
$u^{n+1/3}_1 = F( u^{n+1/3}_2, u^{n+1/3}_3, u^{n+1/3}_4, u^{n+1/3}_5, J^{n+1/3}_s)$
to
$u^{n+1/3}_5 = F( u^{n+1/3}_1, u^{n+1/3}_2, u^{n+1/3}_3, u^{n+1/3}_4, J^{n+1/3}_s)$
which are evaluated after solving a decoupled, closed one-dimensional equation for $J^{n+1/3}_s$. This nonlinear equation is a polynomial with logarithmic terms (from the entropy expression). It is smooth and regular and easily solved by standard root-finding techniques. This approach needs some thorough preparation in a symbolic programming language, but the results for all attempted shock tube simulations are perfectly smooth.

\section{Results}
\label{sec:results}
\begin{figure}[h]
\includegraphics{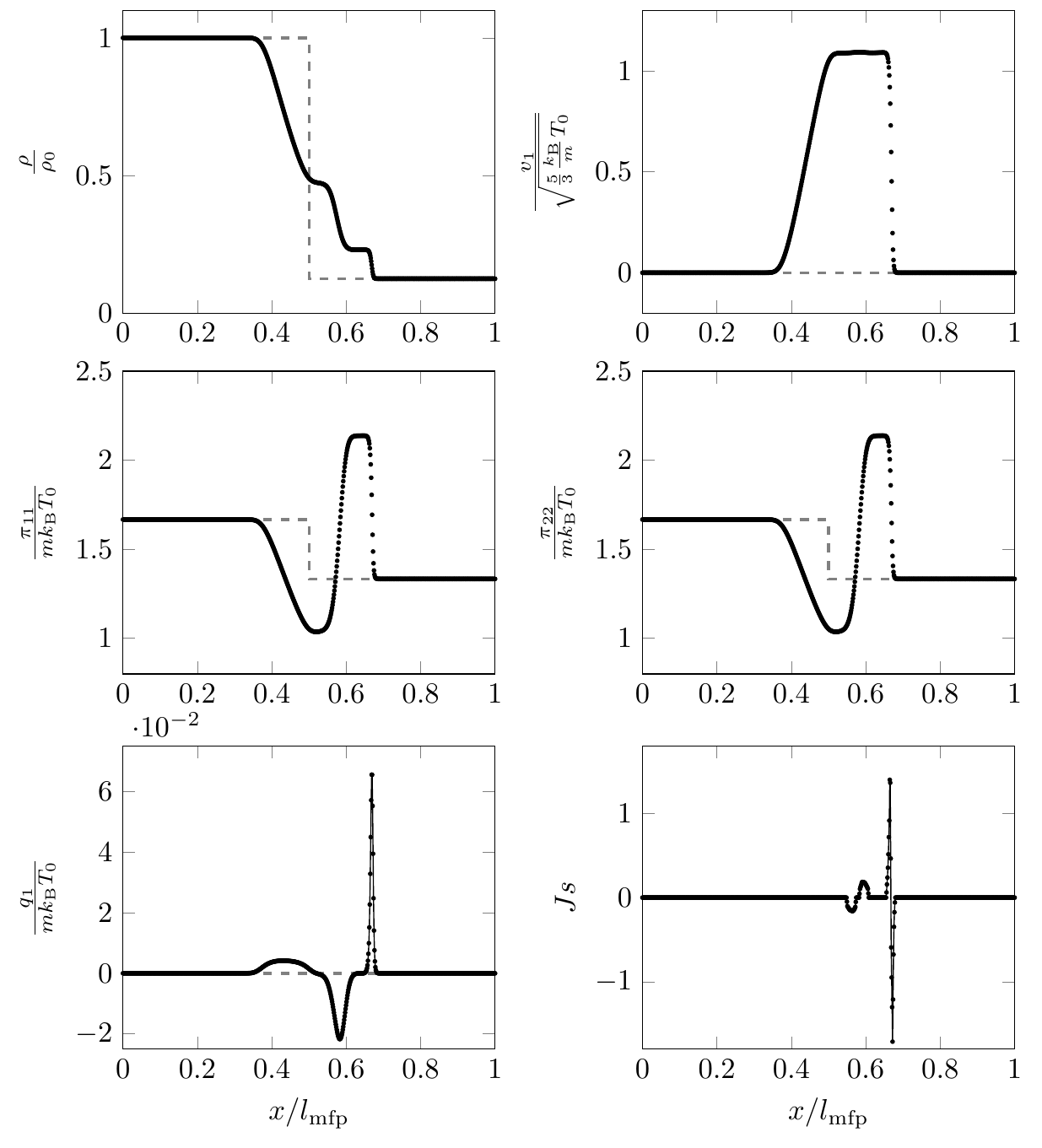}
\caption{Simulation results for $N=800$, $\lambda=1/9$ and $\epsilon=10^{-4}$. The dashed line shows the initial conditions, the dots represent the grid point values of the simulation variables. The continuous line in the last two plots is for visual guidance.}
\label{fig:vars1}
\end{figure}

Since we are applying a novel numerical approach to a new physical model we cannot benchmark the results independently. Therefore we show two simulation results. First we adhere to the simulation setup and parameters proposed in \cite{liotta2000central}, which serves as a general demonstration of the numerical functionality. The initial values for the Riemann problem are
\begin{align}
\rho = 1, \quad \pi_{11}=\pi_{22}= 5/3  \quad&\text{for} \quad0<x<0.5\\
\rho = 1/8, \quad\pi_{11}=\pi_{22}= 4/3 &\quad\text{for} \quad0.5\geq x\geq 1
\end{align}
and $M_1$ and $q_1$ equal to zero. The simulation grid size is $N=800$ and $\lambda=1/9$ is within the stability region. Dimensionless variables are obtained by rescaling the density with a reference upstream density $\rho_0$. The velocity is scaled by the speed of sound $\sqrt{\frac{5}{3} \frac{k_\text{B}}{m}T_0}$ and the space variable by the length of the mean free path $l_\text{mfp}$. The remaining variables follow from this and the resulting rescaling is identical to the mathematical approach given in \cite{au2001shock}. The dimensionless constants  are chosen as $F=5/3$, $b=1/20$ and $\bar{D}=4/3$. The relaxation time $\tau$ is replaced by the relaxation parameter $\epsilon = 10^{-4}$ which represents the Knudsen number and controls the rarefaction regime. The simulation time is $\hat{t} = 0.07$. The results for the simulation variables are presented in Figure ~\ref{fig:vars1} and exhibit perfectly smooth solutions.
At this stage, the values of the new variable $J_s$ have already settled down and only the most prominent peaks remain. Smaller features are no longer visible since the order of $J_s$ in $g(u)$ is $10^{-4}$ smaller than the other production terms in this regime.

\begin{figure}[h]
\includegraphics{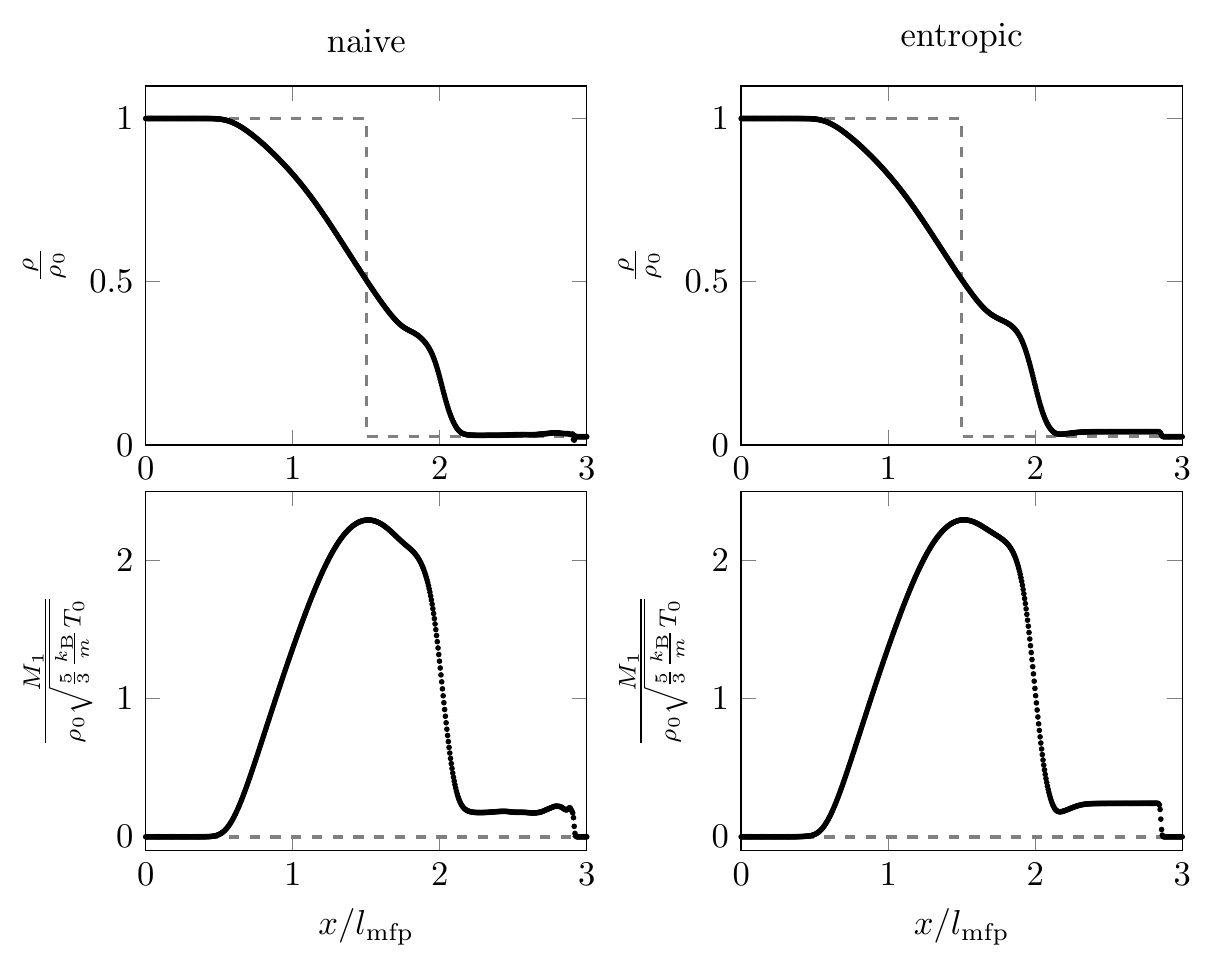}
\caption{Comparison of the naive and entropic schemes for density and momentum. Simulation results for $N=600$, $\lambda=0.025$ and $\epsilon=1$. The dashed line shows the initial conditions, the dots represent the grid point values.}
\label{fig:comparison1}
\end{figure}

\begin{figure}[h]
\includegraphics{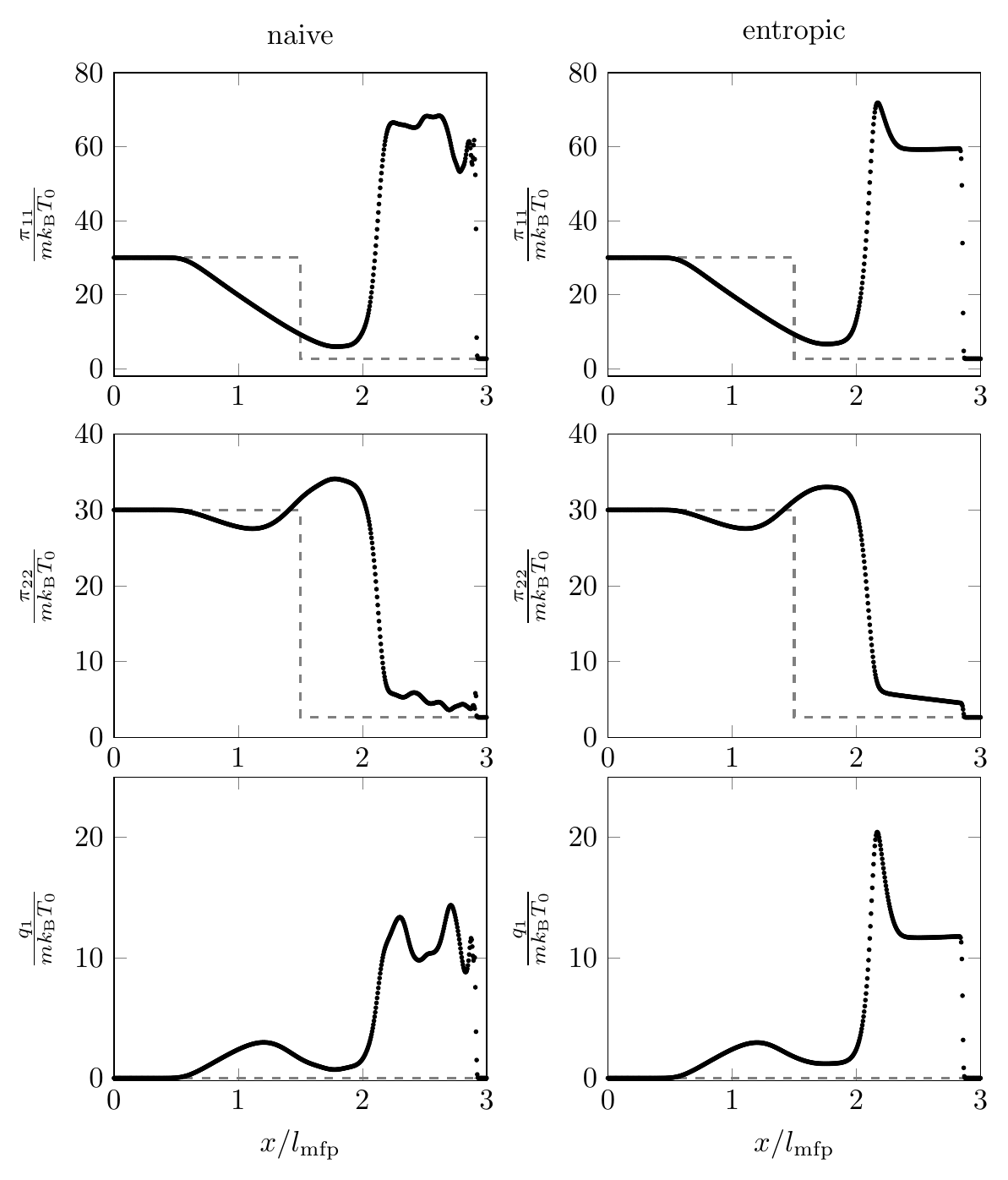}	
\caption{Comparison of the naive and entropic schemes for the second moment tensor and third moment vector. Simulation results for $N=600$, $\lambda=0.025$ and $\epsilon=1$.}\label{fig:comparison2}
\end{figure}

For the comparison of results between the naive and the entropic approach we alter the simulation regime. We have identified two situations for which the entropic approach is clearly superior. First, it can be observed that the naive approach has a smaller region of stability and oscillations arise earlier when the size of the time step is increased. Of greater importance is the second situation. As we increase the difference between the initial conditions of the left and right hand side, say to
\begin{align}
\rho = 1, \quad \pi_{11}=\pi_{22}= 30  \quad&\text{for} \quad0<x<0.5\\
\rho = 1/40, \quad\pi_{11}=\pi_{22}= 8/3 &\quad\text{for} \quad0.5\geq x\geq 1,
\end{align}
we observe higher shock speeds and greater gradients. These cannot be resolved smoothly by the naive approach even for small time steps. For density and momentum density, the oscillations are not as drastic, although it is evident that the shock front is propagated at slightly faster speed (Figure ~\ref{fig:comparison1}). For the remaining variables the difference between the spoilt and smooth simulations are far more striking (Figure ~\ref{fig:comparison2}).


\section{Conclusion}
We have introduced a general concept of using the redundant entropy balance to our advantage for numerical schemes. The method produces perfectly smooth and stable solutions. Although we are demonstrating the principle by showing a particular example, the method is clearly more universal and could be used in other hyperbolic or mixed systems with an entropy evolution equation.

Additionally, the concept could open doors for certain open question in the theory of hyperbolic equations when a balance form is not attainable, such as the development of further schemes and enables the calculation of Rankine-Hugoniot conditions.

\section*{Acknowledgments}
The authors thank Manuel Torrilhon, Henning Struchtrup and Martin Kr\"oger for valuable discussions.


\begin{thebibliography}{10}

\bibitem{au2001shock}
{\sc JD~Au, M~Torrilhon, and W~Weiss}, {\em The shock tube study in extended
  thermodynamics}, Phys.\ Fluids, 13 (2001), p.~2423.

\bibitem{grad}
{\sc Harold Grad}, {\em On the kinetic theory of rarefied gases},
  Communications on Pure and Applied Mathematics, 2 (1949), pp.~331--407.

\bibitem{karlin2006elements}
{\sc Iliya~V Karlin, Santosh Ansumali, Christos~E Frouzakis, and Shyam~Sunder
  Chikatamarla}, {\em Elements of the lattice {B}oltzmann method i: Linear
  advection equation}, Commun.\ Comput.\ Phys, 1 (2006), pp.~616--655.

\bibitem{karlin2007elements}
{\sc Iliya~V Karlin, Shyam~S Chikatamarla, and Santosh Ansumali}, {\em Elements
  of the lattice {B}oltzmann method ii: Kinetics and hydrodynamics in one
  dimension}, Commun.\ Comput.\ Phys, 2 (2007), pp.~196--238.

\bibitem{leveque1992numerical}
{\sc Randall~J LeVeque}, {\em Numerical methods for conservation laws},
  (1992).

\bibitem{liotta2000central}
{\sc Salvatore~Fabio Liotta, Vittorio Romano, and Giovanni Russo}, {\em Central
  schemes for balance laws of relaxation type}, SIAM J.\ Numer.\ Anal.\, 38
  (2000), pp.~1337--1356.

\bibitem{nessyahu1990non}
{\sc Haim Nessyahu and Eitan Tadmor}, {\em Non-oscillatory central differencing
  for hyperbolic conservation laws}, J.\ Comput.\ Phys., 87 (1990),
  pp.~408--463.

\bibitem{ottinger2010ottinger}
{\sc Hans~Christian {\"O}ttinger}, {\em {\"O}ttinger replies}, Phys.\ Rev.\
  Lett., 105 (2010), p.~128902.

\bibitem{paperhco}
\leavevmode\vrule height 2pt depth -1.6pt width 23pt, {\em Thermodynamically
  admissible 13 moment equations from the {B}oltzmann equation}, Phys.\ Rev.\
  Lett., 104 (2010), p.~120601.

\bibitem{hcocomment}
{\sc Henning Struchtrup and Manuel Torrilhon}, {\em Comment on
  thermodynamically admissible 13 moment equations from the {B}oltzmann
  equation}, Phys.\ Rev.\ Lett., 105 (2010), p.~128901.

\end{thebibliography}
\end{document}